\documentclass[11pt, a4paper, UKenglish]{amsart}

\usepackage[utf8]{inputenc}
\usepackage[UKenglish]{babel}
\usepackage{amsmath,amscd,amssymb,amsthm, amsfonts}
\usepackage{amssymb}
\usepackage{latexsym,enumerate,graphicx,geometry,rotating,braket}
\usepackage{mathtools}
\usepackage[all]{xy}
\usepackage[style=alphabetic,natbib=true, backend=biber, bibencoding=utf8]{biblatex}
\usepackage{thmtools}
\usepackage{thm-restate}
\usepackage[colorlinks=true,breaklinks]{hyperref}
\usepackage{cleveref}
\usepackage{tikz-cd}
\usepackage{stmaryrd}
\usepackage{bbm}
\usepackage{mathrsfs}
\usepackage{multirow}
\usepackage{enumitem}

\usepackage{academicons}
\definecolor{orcidlogocol}{HTML}{A6CE39}

\renewcommand{\O}{\mathcal O}
\newcommand{\N}{\mathcal N}

\newcommand{\Q}{\mathbb{Q}}
\newcommand{\Qbar}{\bar{\mathbb{Q}}}

\renewcommand{\t}{\times}

\newcommand{\sL}{\mathcal{L}}
\newcommand{\sM}{\mathcal{M}}

\newcommand{\sO}{\mathcal{O}}

\newcommand{\sX}{\mathcal{X}}

\newcommand{\bbC}{\mathbb{C}}

\newcommand{\bbN}{\mathbb{N}}
\newcommand{\bbP}{\mathbb{P}}
\newcommand{\bbQ}{\mathbb{Q}}
\newcommand{\bbR}{\mathbb{R}}
\newcommand{\bbZ}{\mathbb{Z}}

\newcommand{\into}{\hookrightarrow}

\DeclareMathOperator{\Pic}{Pic}

\DeclareMathOperator{\Spec}{Spec}

\DeclareMathOperator{\End}{End}

\DeclareMathOperator{\M}{M}

\DeclareMathOperator{\an}{an}
\DeclareMathOperator{\ad}{ad}

\DeclareMathOperator{\Disc}{Disc}
\DeclareMathOperator{\Gal}{Gal}

\newcommand{\blank}{{-}}

\DeclarePairedDelimiter\abs{\lvert}{\rvert}%
\DeclarePairedDelimiter\norm{\lVert}{\rVert}%

\makeatletter
\let\oldabs\abs
\def\abs{\@ifstar{\oldabs}{\oldabs*}}
\let\oldnorm\norm
\def\norm{\@ifstar{\oldnorm}{\oldnorm*}}
\makeatother

\newcommand{\defeq}{\vcentcolon=}

\newcommand\restr[2]{{
  \left.\kern-\nulldelimiterspace 
  #1 
  \vphantom{\big|} 
  \right|_{#2} 
  }}
  
\makeatletter
\let\save@mathaccent\mathaccent
\newcommand*\if@single[3]{%
  \setbox0\hbox{${\mathaccent"0362{#1}}^H$}%
  \setbox2\hbox{${\mathaccent"0362{\kern0pt#1}}^H$}%
  \ifdim\ht0=\ht2 #3\else #2\fi
  }
\newcommand*\rel@kern[1]{\kern#1\dimexpr\macc@kerna}
\newcommand*\widebar[1]{\@ifnextchar^{{\wide@bar{#1}{0}}}{\wide@bar{#1}{1}}}
\newcommand*\wide@bar[2]{\if@single{#1}{\wide@bar@{#1}{#2}{1}}{\wide@bar@{#1}{#2}{2}}}
\newcommand*\wide@bar@[3]{%
  \begingroup
  \def\mathaccent##1##2{%
    \let\mathaccent\save@mathaccent
    \if#32 \let\macc@nucleus\first@char \fi
    \setbox\z@\hbox{$\macc@style{\macc@nucleus}_{}$}%
    \setbox\tw@\hbox{$\macc@style{\macc@nucleus}{}_{}$}%
    \dimen@\wd\tw@
    \advance\dimen@-\wd\z@
    \divide\dimen@ 3
    \@tempdima\wd\tw@
    \advance\@tempdima-\scriptspace
    \divide\@tempdima 10
    \advance\dimen@-\@tempdima
    \ifdim\dimen@>\z@ \dimen@0pt\fi
    \rel@kern{0.6}\kern-\dimen@
    \if#31
      \overline{\rel@kern{-0.6}\kern\dimen@\macc@nucleus\rel@kern{0.4}\kern\dimen@}%
      \advance\dimen@0.4\dimexpr\macc@kerna
      \let\final@kern#2%
      \ifdim\dimen@<\z@ \let\final@kern1\fi
      \if\final@kern1 \kern-\dimen@\fi
    \else
      \overline{\rel@kern{-0.6}\kern\dimen@#1}%
    \fi
  }%
  \macc@depth\@ne
  \let\math@bgroup\@empty \let\math@egroup\macc@set@skewchar
  \mathsurround\z@ \frozen@everymath{\mathgroup\macc@group\relax}%
  \macc@set@skewchar\relax
  \let\mathaccentV\macc@nested@a
  \if#31
    \macc@nested@a\relax111{#1}%
  \else
    \def\gobble@till@marker##1\endmarker{}%
    \futurelet\first@char\gobble@till@marker#1\endmarker
    \ifcat\noexpand\first@char A\else
      \def\first@char{}%
    \fi
    \macc@nested@a\relax111{\first@char}%
  \fi
  \endgroup
}

\def\house#1{{%
    \setbox0=\hbox{$#1$}
    \vrule height \dimexpr\ht0+1.4pt width .4pt depth \dp0\relax
    \vrule height \dimexpr\ht0+1.4pt width \dimexpr\wd0+2pt depth \dimexpr-\ht0-1pt\relax
    \llap{$#1$\kern1pt}
    \vrule height \dimexpr\ht0+1.4pt width .4pt depth \dp0\relax
}}

\def\House#1{{%
    \setbox0=\hbox{$#1$}
    \vrule height \dimexpr\ht0+1.4pt width .4pt depth \dp0\relax
    \vrule height \dimexpr\ht0+1.4pt width \dimexpr\wd0+2pt depth \dimexpr-\ht0-1pt\relax
    \llap{$#1$\kern1pt}
    \vrule height \dimexpr\ht0+1.4pt width .4pt depth \dp0\relax
}}

\def\housealp{{%
    \setbox2=\hbox{$\alpha$}
    \vrule height \dimexpr\ht2+1.75pt width .4pt depth \dp2\relax
    \vrule height 6.05pt width \dimexpr\wd2+2pt depth -5.65pt
    \llap{$\alpha$\kern1pt}
    \vrule height \dimexpr\ht2+1.75pt width .4pt depth \dp2\relax
}}

\newcommand{\xdasharrow}[2][->]{
\tikz[baseline=-\the\dimexpr\fontdimen22\textfont2\relax]{
\node[anchor=south,font=\scriptsize, inner ysep=1.5pt,outer xsep=2.2pt](x){#2};
\draw[shorten <=3.4pt,shorten >=3.4pt,dashed,#1](x.south west)--(x.south east);
}
}

\newcommand{\eqn}[1]{\begin{equation*}#1\end{equation*}}

\newtheorem{thm_}{Theorem}[section]
\newtheorem{lemma_}[thm_]{Lemma}
\newtheorem{prop_}[thm_]{Proposition}
\newtheorem{conj_}[thm_]{Conjecture}
\newtheorem{cor_}[thm_]{Corollary}
\theoremstyle{definition}
\newtheorem{eg_}[thm_]{Example}
\newtheorem{def_}[thm_]{Definition}
\newtheorem{rk_}[thm_]{Remark}
\newtheorem{qu_}[thm_]{Question}
\newcommand{\thm}[1]{\begin{thm_}#1\end{thm_}}
\newcommand{\lemm}[1]{\begin{lemma_}#1\end{lemma_}}
\newcommand{\eg}[1]{\begin{eg_}#1\end{eg_}}
\newcommand{\prop}[1]{\begin{prop_}#1\end{prop_}}
\newcommand{\defi}[1]{\begin{def_}#1\end{def_}}
\newcommand{\rk}[1]{\begin{rk_}#1\end{rk_}}
\newcommand{\cor}[1]{\begin{cor_}#1\end{cor_}}

\newcommand{\pf}[1]{\begin{proof}#1\end{proof}}

\title{Fields with few small points}
\author{Nuno Hultberg}


\address{Nuno Hultberg. University of Copenhagen, Institute of Mathematics, Universitetsparken 5, 2100 Copenhagen, Denmark;
ORCiD: \href{https://orcid.org/0000-0003-0097-0499}{orcid.org/0000-0003-0097-0499}
}
\email{nh@math.ku.dk}

\calclayout
\begin{document}

\maketitle

\begin{abstract}
Let $X$ be a projective variety over a number field $K$ endowed with a height function associated to an ample line bundle on $X$. Given an algebraic extension $F$ of $K$ with a sufficiently big Northcott number, we can show that there are finitely many cycles in $X_{\Qbar}$ of bounded degree defined over $F$. Fields $F$ with the required properties were explicitly constructed in \cite{fab} and \cite{oksa}, motivating our investigation. We point out explicit specializations to canonical heights associated to abelian varieties and selfmaps of $\bbP^n$. We apply similar methods to the study of CM-points. As a crucial tool, we introduce a refinement of Northcott's theorem.
\end{abstract}

There have recently been advances on the study of height properties of algebraic extensions of $\Q$ in \cite{fab} and \cite{oksa}. Let $\N$ denote the Northcott number with respect to the logarithmic Weil height. The key result of their work is the following theorem.

\thm{[Theorem 1.3 \cite{oksa}]\label{fields}For every $t \in [0,\infty]$ there exist sequences of prime numbers $(p_i)_{i \in \bbN}$, $(q_i)_{i \in \bbN}$, and $(d_i)_{i \in \bbN}$ such that the field $F = \Q((\frac{p_i}{q_i})^{1/d_i}| i \in \bbN)$ satisfies $\N(F) = t$.
}

\rk{While not stated, everything in \cite{oksa} can be done over an arbitrary number field $K$. For this, think of $K$ as the first step in the tower.}

The full strength of this result is not necessary for our purposes. Instead we opt for the simpler construction of \cite{fab}.

\thm{[Theorem 1.3 \cite{fab}]\label{fab}For every $t \in [0,\infty)$ there exist sequences of prime numbers $(p_i)_{i \in \bbN}$ and $(d_i)_{i \in \bbN}$ such that $p_i^{1/d_i}$ converges to $\exp(2t)$ and the $p_i$ are strictly increasing.

Given such a sequence, the field $F = \Q({p_i}^{1/d_i}| i \in \bbN)$ satisfies $t \leq \N(F) \leq 2t$.
}

We can show the abundance of extensions of $K$ with large Northcott number as a formal consequence of the above theorem, i.e. using it as a blackbox.

\begin{restatable}{rlemm}{fields}
\label{plenty}
Let $C >0$ be a constant and $K$ a number field. Then there exist uncountably many algebraic extensions $F$ of $K$ such that $\N(F) > C$.
\end{restatable}

For fields satisfying the Northcott property the finiteness of cycles of bounded degree and height is known. It is natural to ask whether a similar result can be extended to fields with known Northcott number.

Let $(X,L)$ be a pair consisting of a variety over a number field $K$ and a line bundle on said variety. In order to state our theorems more elegantly, we write $D(V) = (\dim(V) + 1)\deg(V)$ for homogeneous cycles $V$ on $X$. The line bundle implicit in this notation will be clear from context. Going forward, all cycles will be assumed homogeneous and effective throughout the article.

\begin{restatable}{rthm}{main}
\label{main}Let $X$ be a projective scheme over a number field $K$ endowed with an admissible adelically metrized line bundle $\bar{L}$ whose underlying line bundle $L$ is ample. Let $d \in \bbN$ and $C > 0$ be constants. Then there exists a constant $R > 0$ such that, for all algebraic extensions $F$ of $K$, such that its Northcott number satisfies $\N(F) > d(C + R)$, we obtain the following.

There are only finitely many $F$-rational cycles $V$ on $X$ such that $D(V) \leq d$ and $h_{\bar{L}}(V) < CD(V)$.
\end{restatable}

\rk{Regardless of this theorem, we can't expect to have only finitely many subvarieties defined over even a number field $K$ as the Northcott property holds only for subvarieties of bounded degree. An example of the failure of the Northcott property without bound on the degree are the subvarieties $\overline{\{(z,z^n)\}} \subseteq \bbP^2$. They are all distinct, defined over the base field and have canonical height $0$.}

We will now give some specializations of interest with explicit constants.

\begin{restatable}{rthm}{proj}
\label{proj}Consider $\bbP^n$ over a number field $K$ endowed with the canonical toric height $\hat{h}$. Let $d \in \bbN$ and $C > 0$ be constants. Let $F$ be an extension of $K$, such that its Northcott number satisfies 
\eqn{\N(F) > d\left(C + \frac{7}{2} n\log 2  + \sum_{i = 1}^n \frac{1}{2i} +\log 2\right).}
Then there are only finitely many $F$-rational cycles $V$ on $\bbP^n_{K}$ such that $D(V) \leq d$ and $\hat{h}(V) < CD(V)$.
\end{restatable}

\begin{restatable}{rthm}{abvar}
\label{abvar}Let $A$ be an abelian variety of dimension $g$ over a number field $K$ endowed with an ample symmetric line bundle $\sM$. Let $L$ denote the extension of $K$ generated by
\eqn{\ker\left( A \xrightarrow{[16]}A \xrightarrow{p_{\sM}} A^{\vee}\right),}
where $p_{\sM}$ denotes the polarization morphism associated to $\sM$. Then there is an embedding $\Theta$ of $A$ into $\bbP^n$ defined over $L$ with associated line bundle $\sM^{\otimes 16}$. Denote by $h_2$ the $l^2$-logarithmic Weil height and by $\hat{h}_{\sM}$ the canonical height associated to the group structure of $A$.

Let $d \in \bbN$ and $C > 0$ be constants. If $F$ is an extension of $L$, such that its Northcott number satisfies 
\eqn{\N(F) > \frac{d}{16}\left(C + 4^{g+1}h_2(\Theta_{\sM^{\otimes 16}}(0_A)) + 3g \log 2  + \sum_{i = 1}^n \frac{1}{2i}+\log 2 \right),} 
then there are only finitely many $F$-rational cycles $V$ on $A_L$ such that $D(V) \leq d$ and $\hat{h}_{\sM}(V) < CD(V)$. In particular, there are only finitely many torsion points and abelian subvarieties with $D(V) \leq d$ defined over $F$.
\end{restatable}

A similar result may be obtained for dynamical systems on projective space.

\begin{restatable}{rthm}{dyn}
\label{dyn}Let $f: \bbP^n \to \bbP^n$ be a selfmap of degree $D \geq 2$, defined over a number field $K$. Denote by $\hat{h}$ the canonical height associated to $f$ and the tautological line bundle. Let $d \in \bbN$ and $C > 0$ be constants. Let $F$ be an extension of $K$, such that its Northcott number satisfies \eqn{\N(F) > d\left(C +  C_1(n,D)h(f) + C_2(n,D) + \sum_{i = 1}^n \frac{1}{2i}\right),}
where $h(f)$ is the height of the coefficients of $f$ as a projective tuple and
\eqn{C_1(n,D)=5nD^{n+1},\ \ C_2(n,D)=3^n n^{n+1}(2D)^{n2^{n+4}D^n}.}
Then there are only finitely many $F$-rational effective divisors $V$ on $\bbP^n_K$ such that $\deg(V) \leq d$ and $\hat{h}(V) < CD(V)$. In particular, there are only finitely many preperiodic hypersurfaces of degree $\leq d$ defined over $F$.
\end{restatable}

\rk{Based on the ideas in \cite{ingram}, a result that is linear in $\deg(V)$ should be possible in any codimension. At the present moment we may use \cite[Theorem 4.12]{hutz}, which yields a bound exponential in $\deg(V)$.}

\rk{If we restrict to geometrically irreducible closed subsets we can improve the bound on the Northcott number by $d\log 2 $ in Theorems \ref{main}, \ref{proj} and by $d\log 2 /16$ in Theorem \ref{abvar}. The statement of Theorem \ref{dyn} cannot be improved.}

We lastly consider an application to CM points on the modular curve. These are not small points in the usual sense. For this reason it is necessary to consider weighted Weil heights.

\begin{restatable}{rthm}{CM}
\label{CM}
There are uncountably many algebraic field extensions of $\Q$ containing only finitely many CM $j$-invariants.
\end{restatable}

The author is not aware of other examples of infinite algebraic extensions of $\Q$ known to contain only finitely many CM $j$-invariants.

In the first section we introduce Northcott numbers and their behaviour under field extension. Lastly we deduce Lemma \ref{plenty}.

The second section will deal with various notions of height and the bounds on their differences. At the end we will see how Theorems \ref{main} and \ref{proj} follow from these bounds.

The third section contains the applications to abelian varieties and dynamical systems on projective space.

At last, we construct infinite algebraic extensions of $\Q$ over which only finitely many CM points are defined.

\section*{Acknowledgements}

I thank Fabien Pazuki for his guidance and mathematical discussions. I am specially grateful for his suggestion to consider also positive dimensional subvarieties and pointing me to references.

I thank Desir\'ee Gij\'on G\'omez for helpful comments on drafts of this article.

I thank Ricardo Menares for nice conversations at the journ\'ees arithmetiques 2023 and his suggestion to consider CM points.

I lastly thank Martin Widmer for pointing out that \cite{oksa} is more general than I originally credited it to be and the audience at the Atelier ANR J-invariant for feedback on a wrong attribution.

\section{Northcott numbers}

In this section, we introduce Northcott numbers of subsets of $\Qbar$, which allows us to refine Northcott's theorem (see \cite[Theorem 2.1]{northcott}) to a statement on Northcott numbers that we call the Northcott inequality. We conclude the section with a proof of Lemma \ref{plenty}.

\defi{[Northcott number]For a subset $S \subseteq \Qbar$ of the algebraic numbers we define the Northcott number of $S$ with respect to a function $f: \Qbar \to [0,\infty)$ as
\eqn{\N_f(S) = \inf\{t \in [0,\infty) | \# \{\alpha \in S; f(\alpha)< t\} = \infty\}.}
We follow the convention that $\inf \emptyset = \infty$.
We call $\N(S)\in [0,\infty]$ the Northcott number of $S$.
}

\rk{Our main focus is on the case that $f=h$ is the logarithmic Weil height. In this case, we omit the $h$ from the notation.
}

\eg{Let $K$ be a number field. Then by Northcott's theorem $\N(K) = \infty$. On the other hand, $\N(\Qbar) = 0$. }

We now state and prove the Northcott inequality.

\thm{[Northcott inequality]\label{nc} Let $F$ be a field with Northcott number $\mathcal{N}(F) = C$. Then the set of algebraic numbers $X$ of degree $\leq d$ over $F$ satisfies $\mathcal{N}(X)\geq \frac{C-d\log 2}{d2^d}$.}

\pf{Let $\epsilon > 0$. Let $Y_\epsilon$ be the set of algebraic numbers $x$ of height $\leq \frac{C-d\log 2}{d2^d}-\epsilon = B_\epsilon$ satisfying $[F(x):F]\leq d$. It is enough to show that the set $Y_\epsilon$ is finite for any $\epsilon > 0$. Let $x\in Y_\epsilon$. Then the at most $d$ conjugates of $x$ over $F$ are also elements of $Y_\epsilon$. The coefficients of the minimal polynomial of $x$ over $F$ are elementary symmetric functions in these conjugates. We can bound the height of the coefficients by
\eqn{d2^dB_\epsilon + d\log 2  = C - \epsilon d2^d}
using the properties of the height (see \cite[Prop. 1.5.15]{bombgubler}). Let $x, x_1, \dots, x_r \in \Qbar$ and $\sigma \in \Gal(\Qbar/\Q)$, then
\begin{align}
h(\sigma(x)) = h(x) \\
h(x_1+ \dots + x_r) \leq h(x_1) + \dots + h(x_r) + \log r \\ 
h(x_1 \dots x_r) \leq h(x_1) + \dots + h(x_r).
\end{align}
However, by assumption on $F$, there are only finitely many such coefficients, thus showing the finiteness of $Y_\epsilon$.
}

\rk{The optimal bound we may obtain with these methods is $\min_{0 \leq j \leq d} \frac{C-\log{d \choose j}}{{d \choose j}j}$.}

In \cite[Lemma 5]{fab} they notice that the house shares the crucial properties necessary to perform the proof of Theorem \ref{nc}. By combining the ideas of \cite[Lemma 5]{fab} and Theorem \ref{nc} we obtain.

\lemm{Let $f:\Qbar \to [0,\infty)$ be a function. Denote by $\N_f(S)$ the Northcott number of a subset $S \subseteq \Qbar$ with respect to $f$. Suppose that $f$ satisfies
\begin{align}
f(\sigma(x)) = f(x) \\
f(x_1+x_2) \leq F(f(x_1),f(x_2))\\ 
f(x_1 x_2) \leq F(f(x_1),f(x_2))
\end{align}
for some continuous function $F:\bbR^2 \to [0,\infty)$ and all $x_1,x_2 \in \Qbar$ and $\sigma \in \Gal(\Qbar/\Q)$. Then there exists a continuous function $G:[0,\infty] \to [0,\infty]$ with $G(\infty) = \infty$ depending only on $F$ and an auxiliary natural number $d$ such that the following holds. Let $U \subseteq \Qbar$ and let $S\subseteq \Qbar$ be the subset of numbers satisfying monic polynomials with coefficients in $U$ of degree bounded by $d$. Then
\eqn{\N_f(S) \geq G(\N_f(U)).}
}

Let us be more explicit in the case of the house. The house is defined as follows.

\begin{align}
\House{\blank} :\Qbar &\to [0,\infty)\\
\alpha &\mapsto \max_{\sigma:\Qbar \into \bbC}|\sigma(\alpha)|
\end{align}

\lemm{Let $F$ be a field such that $\N_{\house{\blank}}(\O_F) = C$. Then the set of algebraic integers $X$ of degree $\leq d$ over $F$ satisfies $\mathcal{N}_{\house{\blank}}(X)\geq \frac{C^{1/d}}{2^d}$.
}
\pf{The proof is analogous to that of Theorem \ref{nc} using the properties
\begin{align}
\house{\sigma(x)} = \house{x} \\
\house{x_1+x_2} \leq \house{x_1}+ \house{x_2}\\ 
\house{x_1 x_2} \leq \house{x_1}\house{x_2}.
\end{align}
for $x_1, x_2 \in \Qbar$ and $\sigma \in \Gal(\Qbar/\Q)$
}

\rk{We may improve the constant to $\min_{0\leq j\leq d} \frac{C^{1/j}}{{d \choose j}}$.}

This approach, of course, can be used to upper bound Northcott numbers, as well.

\cor{\label{nc2}Suppose a field $K$ has a field extension $F$ of degree $d$ satisfying $\mathcal{N}(F)=C$. Then $\mathcal{N}(K)\leq Cd2^d + d\log 2$.}

\rk{Again we may improve the bound. Here the best possible bound is $\min_{0\leq j\leq d} {d \choose j}j C + \log{d \choose j}$.}

\eg{We may apply this to the field extension $\Q^{tr}(i)/\Q^{tr}$ of the totally real numbers. In \cite[Example 5.3]{bog} it is shown that
\eqn{\alpha_k = \left(\frac{2-i}{2+i}\right)^{1/k}}
is a sequence of points with height tending to zero in $\Q^{tr}(i)$. In particular, $\mathcal{N}_h(\Q^{tr}(i)) = 0$. Hence $\mathcal{N}(\Q^{tr}) \leq \log 2 \approx 0.693$. The best known bound is the one in \cite{tr} ($\mathcal{N}(\Q^{tr}) \leq 0.2732\dots$).}

\rk{The bound in the specific case of the totally real numbers is not sharp and may be improved. Using that the conjugates of $\alpha_k$ equidistribute around the unit circle we may see that $h(\alpha_k + \widebar{\alpha_k}) \to \int_0^1 \max\{2\log\abs{\cos(\pi x)}, 0\} \approx 0.323$.\footnote{This constant also appears as the Mahler measure of the polynomial $1+x+y$,computed by Smyth in \cite{boyd} and as the Arakelov-Zhang pairing $\langle x^2, 1-(1-x)^2\rangle$ in \cite{dynpair}. It equals $\frac{3\sqrt{3}}{4\pi}L(2,\chi)$, where $\chi$ is the nontrivial quadratic character modulo $3$.} }

We can prove lemma \ref{plenty}.

\fields*

\pf{When the ground field is $\Q$, this follows immediately by the work of \cite{oksa} or \cite{fab} quoted at the beginning of the introduction.

Consider now the case of an arbitrary number field $K$ and write $d = [K:\Q]$. We may use Theorem \ref{nc} to obtain that for fields $F$ satisfying $\N(F) > d2^dD+d\log 2$ the composite field $KF$ satisfies $\N(KF) > D$. Over $\Q$, there are uncountably many fields satisfying $\N(F) > d2^dD+d\log 2$. Hence it suffices to show that $KF$ are distinct for distinct $F$.

For this let us consider fields of the form $F = \Q({p_i}^{1/d_i}| i \in \bbN)$, where all $p_i$ and $d_i$ are distinct primes. We can find an extension $F$ of the above form that further satisfies that ${p_i}^{1/d_i}$ tends to $\exp{2t}$ for some $t > d2^dD+d\log 2$. This satisfies the conditions of \ref{fab} and hence $\N(F) \geq t$. Let $t^\prime \neq t$ and $F^\prime$ be an extension $\Q({p^\prime_i}^{1/d^\prime_i}| i \in \bbN)$ with the same conditions of $F$, but with ${p^\prime_i}^{1/d^\prime_i}$ going to $\exp{2t^\prime}$. We need to show that $KF$ cannot contain $F^\prime$. Now $F^\prime$ contains infinitely many ${p^\prime_i}^{1/d^\prime_i}$ that are not contained in $F$. When $d^\prime_i > [K:\Q]$, then also ${p^\prime_i}^{1/d^\prime_i} \notin KF$.
}

\thm{Let $C >0$ be a constant and $K$ a number field. Then there exist uncountably many algebraic extensions $F$ of $K$ such that $\N_{\house{\blank}}(\O_F) > C$.
}

\pf{Fields $F$ with prescribed value for $\N_{\house{\blank}}(\O_F)$ are constructed in \cite[Theorem 1]{fab}. The same argument as above applies since the fields are of similar form.}

\subsection{Relative Northcott numbers}

In \cite{okrelative}, Northcott numbers are considered in a relative setting. The following simplified statement of their result suffices for our needs.

\thm{[\cite{okrelative} Thm. 1.7.]There exists a field $L$ satisfying $\N(L)=0$ such that, for every $t\in(0,\infty]$, there exist sequences of prime numbers $(p_i)_{i \in \bbN}$, $(q_i)_{i \in \bbN}$, and $(d_i)_{i \in \bbN}$ such that the field $F = L((\frac{p_i}{q_i})^{1/d_i}| i \in \bbN)$ satisfies $\N(F\setminus L) = t$.}

\lemm{Let $L \subseteq F \subseteq \Qbar$ be fields satisfying $\N(L)=c$ and $\N(F\setminus L) = t$. Then there exists no $x\in F\setminus L$ satisfying $h(x) < t-c$.}

\pf{We notice that the set $F\setminus L$ is closed under multiplication by elements in $L^\t$. Suppose $x \in F\setminus L$ satisfies $h(x) < t-c$. Let $\epsilon > 0$ be such that $h(x)+2\epsilon < t-c$. Then for any of the infinitely many $y\in L^\t$ satisfying $h(y) \leq c+\epsilon$, $yx$ lies in $F\setminus L$ and satisfies $h(yx) \leq h(y) + h(x) < t-c-2\epsilon+c+\epsilon = t-\epsilon$. This contradicts the assumption $\N(F\setminus L) = t$.}

Using the lemma above we can state and prove our results in a relative setting. Theorem \ref{proj}, for instance, would take the following form.

\thm{Consider $\bbP^n$ over an algebraic extension $L/\Q$ endowed with the canonical toric height $\hat{h}$. Let $d \in \bbN$ and $C > 0$ be constants. Suppose that $\N(L) = c$. Let $F$ be an extension of $K$, such that its relative Northcott number satisfies 
\eqn{\N(F\setminus L) > d\left(C + \frac{7}{2} n\log 2  + \sum_{i = 1}^n \frac{1}{2i} +\log 2\right) + c.}
Then all $F$-rational cycles $V$ on $\bbP^n_{K}$ such that $D(V) \leq d$ and $\hat{h}(V) < CD(V)$ are already defined over $K$.
}

\section{Heights}

This section will contain an overview of some different notions of heights and the bounds on their differences. The two notions of heights we will consider are Arakelov heights, which are defined using arithmetic intersection theory, and Philippon heights, whose definition relies on Chow forms of subvarieties of projective space. While Arakelov heights have conceptual advantages, Philippon height will be crucial to obtain information on the height of a subvariety from the arithmetic of its field of definition.

As a link between these two notions we use canonical heights. Canonical heights may be considered as Arakelov heights, but can at the same time be obtained from Philippon heights by a limit procedure. We will lastly apply this study to prove Theorems \ref{main} and \ref{proj}.

\subsection{Arakelov heights and adelic metrics}

We now introduce the notions in Arakelov geometry needed in this text. For a more comprehensive survey, we refer to \cite{cl}.

Let $X$ be a proper scheme over $\Q$. For all places $v \leq \infty$ we may associate an analytic space $X^{\an}_v$. For $v = \infty$ we set $X^{\an}_\infty = X(\bbC)/F_\infty$, where $F_\infty$ denotes complex conjugation. For $v < \infty$ the definition of the analytification is due to Berkovich in \cite{berk}. For all $v$ this is a compact metrizable, locally contractible topological space containing $X(\bbC_v)/\Gal(\bbC_v/\bbQ_v)$ as a dense subspace. Further, it's equipped with the structure of a locally ringed space with a valued structure sheaf $\sO_{X^{\an}_v}$, i.e. to each $f \in \sO_{X^{\an}_v}(U)$ we can associate an absolute value function $|f|: U \to \bbR_{+}$ that is continuous in a way that is compatible with restrictions. We define $X_{\ad} = \coprod_{v\leq \infty} X^{\an}_v$.

We now define the structure of an adelic metric on a line bundle $L$ on $X$. An \emph{adelic metric} is a collection of compatible $v$-adic metrics. A $v$-adic metric on a line bundle $L^{\an}_v$ on $X^{\an}_v$ is the association of a norm function $||s||_v:U \to \bbR_+$ to every section $s \in L^{\an}_v(U)$ compatible with restriction. Being a norm function means compatibility with multiplication by holomorphic functions and that $||s||_v$ only vanishes when $s$ does. Tensor products and inverses of line bundles with $v$-adic metrics are canonically endowed with $v$-adic metrics. The absolute value endows the trivial bundle with a $v$-adic metric at all places.

The compatibility conditions for adelic metrics reflect the global nature of $X$. A model $(\sX,\sL)$ of $(X,L)$ over $\Spec \bbZ$ induces $v$-adic metrics at all finite places. For a collection of $v$-adic metrics to form an adelic metric we demand it agrees with the metrics induced by $(\sX,\sL)$ at all but finitely many places. If for some power the metrics agree at all places with model metrics we say that the adelic metrics are \emph{algebraic}.

Not all adelically metrized line bundles can be studied equally well. It is often helpful to impose algebraicity and positivity conditions. A notion fulfilling these requirements is semipositivity. \emph{Semipositive metrics} are limits of algebraic metrics with a positivity condition. Important examples of semipositive metrics are the canonical metrics obtained from polarized dynamical systems. An adelic line bundle is called \emph{admissible} if it can be represented as the difference of semipositive adelic line bundles.

We can easily define the height of a point $P \in X(\Qbar)$ in terms of adelic metrics. Let $\bar{L}$ be an adelically metrized line bundle on $X$ with underlying line bundle $L$ and $P \in X(\Qbar)$. This point defines a point $P_v$ in the Berkovich space $X^{\an}_v$ for all $v$. The height of a point $P \in X(\Qbar)$ with respect to an adelically metrized line bundle $\bar{L}$ on $X$ is defined as $h_{\bar{L}}(P) = -\sum_{v \leq \infty} \log||s(P_v)||_v$, where $s$ is a meromorphic section of $L$ with no poles or zeroes at $P$.

More generally, the height of irreducible closed subsets of $X_{\Qbar}$ is defined using arithmetic intersection theory. Given an irreducible closed subset $Z \subseteq X_{\Qbar}$ of dimension $d$, we define
\eqn{h_{\bar{L}}(Z) = \widehat{\deg}_{\bar{L}}(Z)= \widehat{\deg}(\hat{c}_1(\bar{L})^{d+1}|Z).}
We do not follow the convention of \cite{cl} since we would like a notion which is additive in cycles. Our convention differs from that of \cite{cl} by the factor $D(V)= (\dim(V) + 1)\deg(V)$.

\subsection{Heights under the variations of metrics}

We will now introduce a lemma comparing the heights with respect to two admissible metrics.

\lemm{\label{arcomp}Let $X$ be a proper scheme over $\Q$ endowed with a line bundle $L$. Let $\bar{L}$ and $\bar{L}^\prime$ be admissible adelic metrics on $L$. Then there exists a constant $C \in \bbR$ such that for all closed integral subschemes $V \subseteq X_{\Qbar}$ we have
\eqn{|h_{\bar{L}}(V)-h_{\bar{L}^{\prime}}(V)| \leq CD(V).} 
If $L$ is ample and the metrics are algebraic, the admissibility assumption can be omitted.}

\pf{This follows from \cite[Prop. 5.3.]{cl}, a limit argument and linearity. The second case is Prop. 3.7 loc.cit. . In order to follow our convention we multiply the bounds by $D(V)$.}

\subsection{Philippon height}

There is an alternative definition of heights of subvarieties of projective space introduced by Philippon in his papers \cite{phil1}, \cite{phil2} and \cite{phil3}. The Philippon height is obtained from the coefficients of the Chow form of the variety. This viewpoint is important in order to obtain information on the height of a subvariety from the arithmetic of its field of definition. We do not consider the case of weighted projective spaces. For more details we refer to Philippon's original papers. The heights in his different papers differ in the contribution of the infinite places. We will follow \cite{phil3}.

In order to define the Philippon height of a subvariety of projective space we need to first define its Chow form. This is done using projective duality. Let $K$ be a field and $V$ be a closed geometrically irreducible subvariety of $\bbP_K^n$ of dimension $r$. Denote the variety parametrizing linear hyperplanes in $\bbP^n$, i.e. the projective dual of $\bbP^n$, by $\bbP^{n,\vee}$. The subvariety $X$ of $(\bbP^{n,\vee})^{r+1}$ consisting of the tuples of hyperplanes $(H_0, \dots, H_r)$ such that $H_0 \cap \dots H_r \cap V \neq \emptyset$ is a hypersurface. In fact, it is the vanishing locus of a multihomogeneous polynomial over $K$ of degree $\deg V$ in the coordinates of each factor. This polynomial $f$, defined up to multiplication by a scalar, is called the \emph{Chow form} of $V$. If $K$ is a number field we may now proceed to define the Philippon height of $V$. Given the Chow form we define
\eqn{h_{Ph}(V) \defeq \frac{1}{[K:\Q]} \sum_v [K_v:\Q_v] \log{\M_v(f)}.}
Here $\M_v(f)$ is defined as the maximum $v$-adic absolute value of the coefficients of $f$ when $v$ is a finite place. For the archimedean places we define
\eqn{\log\M_v(f) = \int_{(S^{n+1})^{r+1}} \log|\sigma_v(f)| \sigma^{\wedge (r+1)}_{n+1} + D(V)\sum^{n}_{i=1} \frac{1}{2i}.}
Here $\sigma_v$ denotes a choice of complex embedding for the place $v$. $S^{n+1}$ denotes the unit sphere in $\bbC^{n+1}$, while $\sigma_{n+1}$ denotes the invariant probability measure on $S^{n+1}$. We define a variant of the Philippon height $\tilde{h}_{Ph}$ by taking the contribution at an archimedean place to be the maximum modulus of the coefficients instead.

We need to compare the Philippon height with this variant in order to deduce from the Northcott number of a field something about the height of projective varieties defined over said field. Philippon attributes such a comparison to Lelong \cite[Th\'eor\`eme 4]{lelong}. We state it now.

\lemm{\label{comp2} Let $V \subseteq \bbP_{\Qbar}^n$ be an integral closed subvariety, then we have the inequalities
\eqn{0 \leq h_{Ph}(V)-\tilde{h}_{Ph}(V) \leq D(V)\sum_{i=1}^{n}\frac{1}{2i} = D(V)c(n).}
}

Lastly we need to compare Philippon's heights with the toric canonical height on projective space. This allows us to relate Arakelov heights with Philippon heights. The following statement is taken from \cite[Prop 2.1]{philtor}.

\prop{\label{torbound}Let $V \subseteq \bbP_{\Qbar}^n$ be a closed irreducible subset. Let $\hat{h}$ denote canonical toric height on $\bbP^n$. Then
\eqn{|\hat{h}(V)-h_{Ph}(V)| \leq D(V)\frac{7}{2}n \log 2 .}
}

\subsection{\label{cycles}Cycles}

It may be useful to consider the height of general homogeneous cycles defined over a field $F\subseteq \Qbar$. Since the components of an $F$-rational cycle $C$ are not necessarily defined over $F$, a further lemma is required to relate its height to the arithmetic of $F$.

Let $C = \sum n_i V_i$, for geometrically irreducible $V_i$, be a $F$-raitional cycle on $\bbP^n$. Its Chow form is defined to be \eqn{ f_C = \prod f^{n_i}_{V_i}.}
Up to scalar, $f_C$ has coefficients in $F$. Let us define the Philippon height of a cycle $C$ by applying Philippon's construction to $f_C$. We can define $\tilde{h}_{Ph}$ in the analogous way.

The resulting height isn't linear with respect to addition of cycles. To address this issue we invoke an inequality on the height of products of polynomials.

\thm{[\cite{bombgubler} Thm 1.6.13] Let $f_1,\dots,f_m$ be polynomials in $n$ variables, $d$ the sum of partial degrees of $f = f_1\dots f_m$ and let $h$ denote the logarithmic Weil height of the coefficients of a polynomial considered as a projective tuple. Then
\eqn{|h(f)-\sum^m_{j=1} h(f_j)| \leq d\log{2}.}
}

\lemm{\label{cycles}Let $C=\sum n_i V_i$ be a homogeneous cycle of $\bbP_{\Qbar}^n$. Then
\eqn{|\tilde{h}_{Ph}(C) - \sum n_i \tilde{h}_{Ph}(V_i)| \leq D(C)\log 2 .}
}

\pf{We apply the theorem to $f_C = \prod f^{n_i}_{V_i}$ and obtain that $d = (\dim(C) +1)\deg(C)$.}

\subsection{Small subvarieties of projective space}
In this section we prove Theorems \ref{proj} and \ref{main} on small subvarieties.

\proj*

\pf{Let $V=\sum n_i V_i$ be an $F$-rational homogeneous cycle. Then its Chow form $f_V$ has coefficients in $F$. As such, we know that $h(f_V) \leq \N(F)-\epsilon$ for only finitely many cycles. By Lemma \ref{cycles} there can only be finitely many cycles satisfying $\sum n_i \tilde{h}_{Ph}(V_i) \leq \N(F)- D(V)\log 2 -\epsilon$. Consequently there are only finitely many $V$ with $\sum n_i h_{Ph}(V_i)+ \epsilon \leq \N(F) - D(V)\left(c(n) + \log 2 \right)$ by Lemma \ref{comp2}. Moreover, there are only finitely many $V$ such that 
\eqn{\hat{h}(V)+ \epsilon= \sum n_i \hat{h}(V_i)+ \epsilon \leq \N(F) - D(V)\left(\frac{7}{2}n \log 2  +c(n) + \log 2 \right)} by Proposition \ref{torbound}. Under the assumption that $C > \frac{\N(F)}{d} -\frac{7}{2}n \log 2 - c(n)-\log 2$ we obtain that there are only finitely many $F$-rational cycles $V$ on $\bbP^n_K$ such that $D(V) \leq d$ and $\hat{h}(V) < CD(V)$. By rearranging the inequality, we conclude the theorem.
}

We easily obtain Theorem \ref{main} as a consequence.

\main*

\pf{We need to compare the heights on $X$ with heights on projective varieties. For this we replace $\bar{\sL}$ by its $n$-th power such that the underlying line bundle is very ample. Let $X \into \bbP^k$ be an embedding associated to $\sL$. Pulling back the canonical toric metric on $\O(1)$ induces an adelic metric on $\sL$, which we denote $\tilde{\sL}$.

Then by Lemma \ref{arcomp} the height associated to $\tilde{\sL}$ only differs from the one associated to $\bar{\sL}$ by an amount bounded by $R^\prime D(V)$ for some constant $R^\prime$. Now the result follows from Theorem \ref{proj}.}

\rk{As an alternative to admissiblity one may require algebraicity in the above theorem.}

\section{Applications to dynamical systems}

Specializations of our main theorem can be obtained by applying more specific height bounds. The arguments required to obtain these specializations are adaptations of the proof of theorem \ref{proj}, which will only be sketched.

The dynamical systems to be considered in greater detail are the ones given by multiplication on abelian varieties and selfmaps of projective space. We start out with a more general situation considered in the foundational paper of Call and Silverman(\cite{canonical}).

In their setup, $X$ is a smooth projective variety over a number field $K$ endowed with a selfmap $\phi$ and a divisor class $\eta \in \Pic(X) \otimes \bbR$ satisfying $\phi^* \eta = \alpha \eta$ for some $\alpha > 1$. Suppose $h$ is a Weil function associated with $\eta.$ Then there is a constant $R$ such that $|h \circ \phi - \alpha h| \leq R$. Let $\hat{h}$ denote the canonical height for $\eta$ and $\phi$. Then the following holds.

\prop{\label{canonical}[\cite{canonical}Proposition 1.2]For every $P \in X(\bar{K})$, the following inequality holds:
\eqn{|\hat{h}(P) - h(P)| \leq \frac{R}{\alpha - 1}.}
}

Note that we can't expect to have finitely many small points for arbitrary $\eta$, as an associated Weil function might not even be bounded below. We may, however, by adapting the proof of Theorem \ref{proj} obtain the following statement.

\prop{In the current setting, suppose that $\eta$ is very ample and $h$ is induced by the canonical toric height under some embedding into projective space. Let $F$ be an algebraic extension of $K$ satisfying $\N(F) > C + \frac{R}{\alpha - 1}$. Then there are only finitely many points $P \in V(F)$ such that $\hat{h}(P) \leq C$.
}

\pf{We adapt the proof of Theorem \ref{proj}. We bound the height of a point in projective space from below by the height of one of its coordinates and use the bound in Proposition \ref{canonical}.}

\subsection{Small subvarieties of abelian varieties}

In order to study small points on abelian varieties, we embed them into projective space using a variant of the theta embedding, first introduced in \cite{mumford}. For a more detailed overview of its properties, see \cite{subvar}. We will then apply a bound on the difference of the canonical height to the Philippon height from loc.cit. to deduce a result on small points of abelian varieties.

Let $A$ be a $g$-dimensional abelian variety defined over a number field $K$. Let $\sM$ be an ample symmetric line bundle on $A$. Then $\sM^{\otimes 16}$ is very ample. David and Philippon choose sections that yield the embedding $\Theta_{\sM^{\otimes 16}}$, or simply $\Theta$, into $\bbP^N$. It is inspired by the embedding of Mumford in \cite{mumford}, but differs from it. As such, it is not defined over $K$ itself, but over the field generated by 
\eqn{\ker\left( A \xrightarrow{[16]}A \xrightarrow{p_{\sM}} A^{\vee}\right),}
where $p_{\sM}$ denotes the polarization morphism associated to $\sM$.

In this setting, we have the following comparison of heights.

\prop{[\cite{subvar}Proposition 3.9.]\label{comp1}Let $V$ be an integral closed subvariety of $A_{\bar{K}}$ and let $h_2$ denote the $l^2$-logarithmic Weil height. Then
\eqn{|\hat{h}_{\sM^{\otimes 16}}(V) - h_{Ph}(\Theta(V))| \leq c_0(\Theta)D(V).}
Here, $c_0(\Theta) = 4^{g+1}h_2(\Theta(0_A)) + 3g \log 2 $.
}

\abvar*

\pf{We adapt the proof of Theorem \ref{proj}. The main differences are that Proposition \ref{comp1} applies to $\hat{h}_{\sM^{\otimes 16}} = 16 \hat{h}_{\sM}$ instead of directly to $\hat{h}_{\sM}$ and that the $\Theta$-embedding of $A$ is not defined over its field of definition $K$, but only over $L = K\left(\ker\left( A \xrightarrow{[16]}A \xrightarrow{p_{\sM}} A^{\vee}\right)\right)$.}

\rk{The $l^2$-logarithmic Weil height $h_2(\Theta_{\sM^{\otimes 16}}(0_A))$ in the theorem is compared to the Faltings height of the abelian variety in \cite{philfal}. This allows for a phrasing of the theorem that does not reference the theta embedding. In \cite{subvar} the quantity $h(\Theta_{\sM^{\otimes 16}}(0_A))$ is denoted by $h(A)$ which may lead to confusion with the Philippon height of $A$, see \cite[Notation 3.2.]{subvar}.}

\subsection{Small subvarieties with respect to dynamical systems on $\bbP^n$}

Another case in which explicit bounds on difference of heights exist are divisors on $\bbP^n$ with a canonical height from a selfmap. In fact, \cite{ingram} proves the following statement.
\thm{\label{compdyn}Let $f: \bbP^n \to \bbP^n$ be a morphism of degree $d \geq 2$ defined over $\Qbar$. Let $V$ be an effective divisor on $\bbP^n$, then
\eqn{|\hat{h}_f(V) - h_{Ph}(V)| \leq (C_1(n,d)h(f) + C_2(n,d))D(V),
}
where $h(f)$ is the height of the coefficients of $f$ as a projective tuple. Moreover, one may choose
\eqn{C_1(n,d)=5nd^{n+1},\ \ C_2(n,d)=3^n n^{n+1}(2d)^{n2^{n+4}d^n}.}
}
For simplicity he states the theorem only for hypersurfaces, but claims there to be no conceptual obstruction to its generalization.

This leads to Theorem \ref{dyn}.

\dyn*

\pf{We adapt the proof of Theorem \ref{proj}. Note that Theorem \ref{compdyn} applies directly to cycles, so the results in section \ref{cycles} are not needed.}

\section{Application to special points}

While special points on Shimura varieties are not small in the usual sense, our approach can still deduce a finiteness result for CM points on the modular curve defined over certain infinite extensions. To this end, we will use weighted Weil heights.

We have some information on the height of special points on the modular curve from \cite{specialpts}. The result on the degree is a restating of the Brauer-Siegel theorem.

\prop{[Proposition 2.1 \cite{specialpts}] \label{spts} Let $x\in \Qbar$. If the elliptic curve $E_x$ of $j$-invariant $x$ has complex multiplication we denote $\Delta(x) = |\Disc(\End(E_x))|$.
\begin{enumerate}
 \item If x is CM, then $[\Q(x):\Q] = \Delta(x)^{1/2 +o(1)}$.
 \item\label{th:second} There exists an effectively computable constant $C$ such that if $x$ is CM, $h(x) \leq \pi \Delta(x)^{1/2}+ C$.
\end{enumerate}
}

\rk{\label{house_enough}In fact, the proof of part \ref{th:second} computes the asymptotic of the house as the discriminant grows: $\house{x} \approx \exp(\pi \Delta(x)^{1/2})$.}

Let $\gamma \in \bbR$, $x\in \Qbar$. Then, the weighted Weil height $h_{\gamma}$ is defined by $h_\gamma(x) = \deg(x)^{\gamma}h(x)$. We may consider Northcott numbers of subsets $S \subseteq \Qbar$ for varying $\gamma$. For a set $S \subseteq \Qbar$, define the sets
\eqn{I_0(S) = \{\gamma\ |\ \N_{h_\gamma}(S) = 0\}, \ \ I_\infty(S) = \{\gamma\ | \ \N_{h_\gamma}(S) = \infty\}.}

We can summarize the work of \cite{oksa} as follows.

\thm{\label{wt}The sets $I_0(S)$ and $I_\infty(S)$ are (in the case of $I_0(S)$ possibly empty) rays. They satisfy $(1, \infty)\subseteq I_\infty$ and $I(S)=\sup I_0 =\inf I_\infty$. For $\gamma \in (-\infty,1)$ and $c \in [0,\infty]$ one can construct a field $F$ such that $I(F) = \gamma$ and $\N_{h_\gamma}(F) = c$.}

We phrase a Corollary of Theorem \ref{spts} in terms of weighted Weil heights.

\cor{The set of CM points $S$ satisfies $I(S) \leq -1$.}

\CM*

\pf{Any field $F$ satisfying $I(F)<1$ constructed in Theorem \ref{wt} fulfills the conditions.}

\rk{Using Remark \ref{house_enough} we see that the corresponding properties for the weighted house  suffice for the conclusion. However, the counterpart to Theorem \ref{wt} has not yet been proven in this setting.}

\medskip
\sloppy
\emergencystretch=1em
\printbibliography

\end{document}